\documentclass[11pt]{article}

\usepackage[utf8]{inputenc}
\usepackage{amsmath}
\usepackage[english]{babel}
\usepackage{amssymb}
\usepackage{amsthm}
\usepackage{mathtools}
\usepackage{esdiff}
\usepackage{microtype}
\usepackage{paralist}
\usepackage{skak}
\usepackage{bbm}
\usepackage[paper=a4paper,left=25mm,right=25mm,top=25mm,bottom=25mm]{geometry}
\usepackage{stmaryrd}
\usepackage{hyperref}
\usepackage{color} 
\usepackage{fancyhdr}
\usepackage{nomencl}
\usepackage{cite}
\usepackage{longtable}
\usepackage{thmtools}
\usepackage{cuted}
\usepackage{MnSymbol}
\declaretheoremstyle[headfont=\normalfont\bfseries]{bfthmstyle}
\declaretheorem[numberwithin=section,style=bfthmstyle]{Theorem}
\declaretheorem[sharenumber=Theorem,style=bfthmstyle]{Lemma}

\newtheorem{prop}[Theorem]{Proposition} 

\declaretheoremstyle[headfont=\normalfont\bfseries,qed={$\diamond$}]{bfdefstyle}
\declaretheorem[sharenumber=Theorem,style=bfdefstyle]{Definition}
\declaretheorem[sharenumber=Theorem,style=bfdefstyle]{Definition-Proposition}
\declaretheorem[sharenumber=Theorem,style=bfdefstyle]{Definition-Lemma}
\declaretheorem[sharenumber=Theorem,style=bfdefstyle]{Example}

\newtheorem{definition-proposition}[Theorem]{Definition-Proposition} 
\newtheorem{definition-lemma}[Theorem]{Definition-Lemma} 

\declaretheoremstyle[headfont=\normalfont\bfseries,qed={$\diamond$}]{bfremstyle}
\declaretheorem[sharenumber=Theorem,style=bfdefstyle]{Remark}

\declaretheoremstyle[
notefont=\normalfont, notebraces={}{},
headformat=\mathbb{N}UMBER~\mathbb{N}AME~\mathbb{N}OTE
]{nopar}

\numberwithin{equation}{section}

\makenomenclature

\author{Jonas Kirchhoff, Bernhard Maschke}
\title{On the Generating Functions of Irreversible port-Hamiltonian Systems}
\date{\today}

\begin{document}
\setlength{\parindent}{0em}
\pagestyle{fancy}
\lhead{Jonas Kirchhoff, Bernhard Maschke}
\rhead{Generating Functions}

\maketitle

\paragraph{Abstract} We study the geometric structure of the drift dynamics of Irreversible port-Hamiltonian systems. This drift dynamics is defined with respect to a product of quasi-Poisson brackets, reflecting the interconnection structure and the constitutive relations of the irreversible phenomena occuring in the system. We characterize this product of quasi-Poisson brackets using a covariant 4-tensor and an associated function. We derive various conditions for which this 4-tensor and the associated function may be reduced to a product of quasi-Poisson brackets.\vspace{-1mm}
 
\vspace{-1mm}

\vfill
\par\noindent\rule{5cm}{0.4pt}\\
\begin{footnotesize}
Corresponding author: Jonas Kirchhoff\\[1em]
Jonas Kirchhoff\\
Institut für Mathematik, Technische Universität Ilmenau, Weimarer Stra\ss e 25, 98693 Ilmenau, Germany\\
E-mail: jonas.kirchhoff@tu-ilmenau.de\\[1em]

Bernhard Maschke\\
Univ. Lyon, Université Claude Bernard Lyon 1, CNRS, LAGEPP UMR 5007, France\\
E-mail: bernhard.maschke@univ-lyon1.fr\\[1em]

Jonas Kirchhoff thanks the Technische Universität Ilmenau and the Freistaat Thüringen for their financial support as part of the Thüringer Graduiertenförderung.

\end{footnotesize}

\newpage

\section{Introduction}

Irreversible port-Hamiltonian Systems \cite{MaRaSb:13b} differ intrinsically from other dynamical models of irreversible thermodynamic systems, in the sense that they are defined with respect to a product of quasi-Poisson brackets and not a symmetric non-negative bracket such as in e.g. \cite{BeEd:91a,Morrison86}. However the properties of this product of brackets are poorly understood. In this paper, we shall associate a four-tensor and a corresponding set of functions with this product and analyze their relation.

\section{Irreversible port-Hamiltonian Systems}

Irreversible port-Hamiltonian Systems are an extension of port-Hamiltonian
Systems which has been defined to represent not only the energy conservation
but also the irreversible entropy creation \cite{MaRaSb:13a,MaRaSb:13b}
and is defined as follows.
\begin{Definition}
\label{def:IPHS} \cite{MaRaSb:13b} An Irreversible port-Hamiltonian
System (IPHS) is a nonlinear control system 
\begin{equation}\label{IPHS}
\begin{aligned}
\frac{dx}{dt}=\gamma\left(x,\mathrm{d}H(x)\right)\left\{ S,H\right\} _{J}J\mathrm{d}H(x)+W\left(x,\mathrm{d}H(x)\right)\\+g\left(x,\mathrm{d}H(x)\right)u,
\end{aligned}
\end{equation}
where $x\left(t\right)\in\mathbb{R}^{n}$ is the state vector, $u\left(t\right)\in\mathbb{R}^{m}$
is the control input, and $dH(x)$ denotes the differential of the function $H$ at $x$ in standard coordinates, defined by

(i) two (smooth) real functions called \emph{Hamiltonian function} $H(x)\in C^{\infty}(\mathbb{R}^{n})$
and \emph{entropy function} $S(x)\in C^{\infty}(\mathbb{R}^{n})$, 

(ii) the \emph{structure matrix} $J\in\mathbb{R}^{n\times n}$
which is constant and skew-symmetric, defining the Poisson bracket  
\[
\left\{ S,H\right\} _{J} := \mathrm{d}S^{\top}\left(x\right)J\mathrm{d} H
\]

(iii) a real function $\gamma(x,\mathrm{d}H(x))=\hat{\gamma}(x)\in C^{\infty}(\mathbb{R}^{n})$,
strictly positive function of the states and co-states

(iv) the vector field $W(x,\mathrm{d} H(x))\in\mathbb{R}^{n}$
and matrix field $g(x,\mathrm{d}H(x))\in\mathbb{R}^{n\times m}$
defining the input map. 
\end{Definition}

Note that the drift dynamics, which may be rewritten for any function
$f\in\mathcal{C}^{\infty}(\mathbb{R}^{n})$
\begin{equation}
\frac{df}{dt}=\gamma\left(x,\mathrm{d} H\right)\left\{ S,H\right\} _{J}\left\{ f,H\right\} _{J}\label{AutonomousIPHS}
\end{equation}
 is defined in such a way that the two axioms of Irreversible Thermodynamics
are satisfied for the isolated thermodynamic system
\begin{itemize}
\item energy conservation: $\frac{dH}{dt}=0$
\item irreversible entropy production: 
\[
\frac{dS}{dt}=\gamma\left(x,\mathrm{d} H\right)\left\{ S,H\right\} _{J}^{2}\geq 0
\]
 
\end{itemize}
Note that taking account of the input leads to the energy balance equation
\[
\frac{dH}{dt}=\mathrm{d}H^{\top}\left(W+gu\right)
\]
 and the entropy balance equation 
\begin{equation}
\frac{dS}{dt}=\underbrace{\gamma\left(x,\mathrm{d}H\right)\left\{ S,H\right\} _{J}^{2}}_{=\sigma_{int}\geq0}+\mathrm{d}H^{\top}\left(W+gu\right)\label{eq:EntropyBalEq_IPHS}
\end{equation}

The reader may find various examples of irreversible systems fitting
into this frame and ranging from heat exchangers to chemical reaction
dynamics and the gas-piston system \cite{MaRaSb:13b,MaRaSb:13a}.
In these physical systems, the bracket $\left\{ S,H\right\} _{J}$
appears to be the driving force of the irreversible phenomenon such
as the temperature difference for the heat conduction or the chemical
affinity for the chemical reaction dynamics. The function $\gamma(x,\mathrm{d}H(x))=\hat{\gamma}(x)$,
corresponds to the constitutive relation of the phenomenon such as
Fourier's law or the chemical reaction kinetics.

Examining the drift dynamics, it may be observed that it is defined by
two functions, the energy and the entropy function and is generated
by a product of Poisson brackets. This departs radically from all alternative suggestions
of structured description of irreversible thermodynamical systems
ranging from the single generator function bracket approach in \cite{BeEd:91a,BeEd:94,Morrison86}
to the two generator functions approach in \cite{Ka:84,GrOt:97a,GrOt:97b}
where the irreversible phenomena (and the irreversible entropy creation)
are described by a \emph{symmetric positive bracket} associated with
a (pseudo-)gradient dynamics.

In this paper, we shall discuss the definition of the drift dynamics
and ask the question whether the functions $E_{J}:=(f,S,H)\mapsto\gamma\{S,H\}_{J}\{f,H\}_{J}$
appearing in the definition of the dynamics \eqref{AutonomousIPHS} is
the only expression with minimal homogeneity that depend only on gradients of $f,S$ and $H$
so that $H$ is preserved and $S$ increases (nontrivially) monotonically along trajectories
of the system.

\section{Conservative-irreversible functions}

Recall first the classical definition of a derivation.

\begin{Definition}[{Derivation, see~\cite[p.\,39]{Co:85}}] \label{def:derivation}
A function $\delta:\mathcal{C}^\infty(\mathbb{R}^n)\to\mathcal{C}^\infty(\mathbb{R}^n)$ is a derivation if, and only if, $\delta$ is linear and fulfills the Leibniz rule,
\begin{align*}
\forall f,g\in\mathcal{C}^\infty(\mathbb{R}^n): \delta(fg) = f\delta(g) + \delta(f)g.
\end{align*}
\end{Definition}

And define \emph{conservative-irreversible functions}, as follows.

\begin{Definition}\label{def:cons_irrev}
A function
\begin{align*}
E:\left(\mathcal{C}^\infty(\mathbb{R}^n)\right)^3\to\mathcal{C}^\infty(\mathbb{R}^n)
\end{align*}
is \textit{conservative-irreversible} if, and only if, there is a function
\begin{align*}
e:\left(\mathcal{C}^\infty(\mathbb{R}^n)\right)^4\to\mathcal{C}^\infty(\mathbb{R}^n)
\end{align*}
with the properties
\begin{enumerate}
\item[(i)] $e$ is a derivation in each of the four arguments (a \emph{four-derivation}),
\item[(ii)] $e$ is symmetric in the third and fourth argument,
\item[(iii)] $e(h,\cdot,h,h) \equiv 0$ for all $h\in\mathcal{C}^\infty(\mathbb{R}^n)$,
\item[(iv)] $e(\cdot,\cdot,h,h)$ is pointwise symmetric and positive semidefinite for all $h\in\mathcal{C}^\infty(\mathbb{R}^n)$,
\end{enumerate}
so that 
\begin{align}\label{eq:tensor-relation}
\forall f,s,h\in\mathcal{C}^\infty(\mathbb{R}^n):E(f,s,h) = e(f,s,h,h).
\end{align}
\end{Definition}

\begin{Example}\label{ex:an_example}
Let $J\in\mathcal{C}^\infty(\mathbb{R}^n,\mathbb{R}^{n\times n})$ be pointwise skew-symmetric. Consider the function
\begin{align*}
e_J:\left(\mathcal{C}^\infty(\mathbb{R}^n)\right)^3 & \to\mathcal{C}^\infty(\mathbb{R}^n),\\
 & (f,s,h,g)\mapsto  \{s, g\}_J \{f,h\}_J.
\end{align*}
This function has the properties (i)-(iv) and induces the conservative-irreversible function
\begin{align*}
E_J:\left(\mathcal{C}^\infty(\mathbb{R}^n)\right)^3 & \to\mathcal{C}^\infty(\mathbb{R}^n),\\
 & (f,s,h)\mapsto  \{s, h\}_J \{f,h\}_J.
\end{align*}
\hfill$\diamond$
\end{Example}

\begin{Remark}
The symmetry condition (ii) states that we choose a symmetric representative of all four-derivations that generate a given conservative-irreversible function $E$. Indeed, if $e$ and $e'$ are any four-derivations with the properties (iii) and (iv) which fulfill~\eqref{eq:tensor-relation}, then, for all $\lambda\in\mathbb{R}$, the four-derivation
\begin{align*}
(f,s,h,g)\mapsto \lambda e(f,s,h,g)+ (1-\lambda)e'(f,s,g,h)
\end{align*}
has the properties (iii) and (iv) and fulfills~\eqref{eq:tensor-relation}, too. Especially, the symmetrisation of $e$,
\begin{align*}
(f,s,h,g)\mapsto \frac{1}{2} e(f,s,h,g)+ \frac{1}{2}e(f,s,g,h),
\end{align*}
is a generating four-derivation for $E$.
\end{Remark}

Algebraically, the set of conservative-irreversible functions has naturally the structure of a module over the smooth, pointwise nonnegative functions, as can easily been seen.

\begin{Lemma}
Equipped with pointwise addition and multiplication, the set of nonnegative smooth functions
\begin{align*}
\mathcal{C}^\infty(\mathbb{R}^n,\mathbb{R}_{\geq 0}):=\left\lbrace f\in\mathcal{C}^\infty(\mathbb{R}^n)\,\big\vert\,\forall x\in\mathbb{R}^n: f(x)\geq 0\right\rbrace.
\end{align*}
is a commutative ring with 1.
\end{Lemma}
\textit{Proof.} Since $\mathbb{R}$ is an ordered field and since the product of smooth real valued functions is smooth (see~\cite[p.\,86]{AbMaRa:88}), $\mathcal{C}^\infty(\mathbb{R}^n,\mathbb{R}_{\geq 0})$ is indeed a commutative ring with $1$.\hfill$\square$

Denote the set of all conservative-irreversible functions with
\begin{align*}
\mathfrak{CI}(\mathbb{R}^n) := \left\lbrace E:\left(\mathcal{C}^\infty(\mathbb{R}^n)\right)^3 \to\mathcal{C}^\infty(\mathbb{R}^n)\,\big\vert\, E~\text{cons.~irrev.}\right\rbrace.
\end{align*}
and equip it with the pointwise addition
\begin{align*}
(E+G)(f,s,h)(x) := E(f,s,h)(x)+G(f,s,h)(x),
\end{align*}
and the multiplication
\begin{align*}
(\lambda\cdot E)(f,s,h)(x) := \lambda(x)E(f,s,h)(x)
\end{align*}
for all $E\in\mathfrak{CI}(\mathbb{R}^n)$ and $\lambda\in\mathbb{C}^\infty(\mathbb{R}^n,\mathbb{R}_{\geq 0})$.

\begin{Lemma}\label{lem:module}
$(\mathfrak{CI}(\mathbb{R}^n),+,\cdot)$ is a module over $\mathcal{C}^\infty(\mathbb{R}^n,\mathbb{R}_{\geq 0})$.
\end{Lemma}
\textit{Proof.}
Let $E,G\in\mathfrak{CI}(\mathbb{R}^n)$ and $\lambda\in\mathcal{C}^\infty(\mathbb{R}^n,\mathbb{R}_{\geq 0})$ be arbitrary. We show that $\lambda E + G\in\mathfrak{CI}(\mathbb{R}^n)$. Since $E$ and $G$ are conservative-irreversible functions, there are functions $e$ and $g$ as in Definition~\ref{def:cons_irrev} so that
\begin{align*}
E = e(f,s,h,h),\quad G(f,s,h) = g(f,s,h,h)
\end{align*}
for all $f,s,h\in\mathcal{C}^\infty(\mathbb{R}^n)$. Then, $\lambda E + G$ has the representation $\lambda e + f$. It remains to verify that $\lambda e + f$ has the properties (i)--(iv) in Definition~\ref{def:cons_irrev}. (i) holds, since derivations on $\mathcal{C}^\infty(\mathbb{R}^n)$ are a module over $\mathcal{C}^\infty(\mathbb{R}^n)$. Let $f,s,h,i\in\mathcal{C}^\infty(\mathbb{R}^n)$. Then
\begin{align*}
(\lambda f + g)(f,s,h,i) & = \lambda f(f,s,h,i) + g(f,s,h,i)\\
& = \lambda f(f,s,i,h) + g(f,s,i,h)\\
& = (\lambda f+g)(f,s,i,h)
\end{align*}
and hence  (ii) holds. The remaining properties can be shown analogously from the pointwise definition of the addition and scalar multiplication on $\mathfrak{CI}(\mathbb{R}^n)$. This shows that the operations $+$ and $\cdot$ are well-defined. It is easy to verify that $(\mathfrak{CI}(\mathbb{R}^n),+,\cdot)$ is indeed a module over $\mathcal{C}^\infty(\mathbb{R}^n,\mathbb{R}_{\geq 0})$; we omit the details.\hfill$\square$

Besides being a module, $\mathfrak{CI}(\mathbb{R}^n)$ is a convex cone.

\begin{Lemma}
$\mathfrak{CI}(\mathbb{R}^n)$ is a convex cone in the real vector space of functions from $\left(\mathcal{C}^\infty(\mathbb{R}^n)\right)^3$ to $\mathcal{C}^\infty(\mathbb{R}^n)$.
\end{Lemma}
\textit{Proof.}
Since constant functions are smooth, this is a direct consequence of Lemma~\ref{lem:module}.\hfill $\square$

We want to characterise the functions $e$ that induce conservative-irreversible functions as in~\eqref{eq:tensor-relation}.

\begin{prop}\label{prop:tensorial_representation}
A four-linear function
\begin{align*}
e:\left(\mathcal{C}^\infty(\mathbb{R}^n)\right)^4\to\mathcal{C}^\infty(\mathbb{R}^n)
\end{align*}
has the properties (i)--(iv) from Definition~\ref{def:cons_irrev} if, and only if, there is some $\varepsilon\in\mathcal{C}^\infty(\mathbb{R}^n,\mathbb{R}^{n\times n\times n\times n})$
with the properties
\begin{enumerate}
\item[(a)] $\varepsilon_{i,j,k,\ell} = \varepsilon_{i,j,\ell,k}$,
\item[(b)] $\varepsilon_{i,j,k,\ell} + \varepsilon_{k,j,\ell,i}+\varepsilon_{\ell,j,i,k} = 0$,
\item[(c)] the function
\begin{align*}
h\mapsto \left[\sum_{k,\ell = 1}^n\varepsilon_{i,j,k,\ell}\mathrm{d}_k h \mathrm{d}_\ell h\right]_{i,j = 1}^n
\end{align*}
is pointwise symmetric and positive semidefinite,
\end{enumerate}
for all $i,j,k,\ell\in\underline{n}$ so that
\begin{align}\label{eq:conservative_tensor}
e(f,s,h,q) = \sum_{i,j,k,\ell = 1}^n \varepsilon_{i,j,k,\ell}\mathrm{d}_i f\mathrm{d}_j s\mathrm{d}_k h\mathrm{d}_\ell q.
\end{align}
\end{prop}
\textit{Proof.} We split the proof into steps.

\textsc{Step 1:} We show that $e$ has a tensor representation~$\varepsilon$. By taking the basis of the tangent space associated to the standard coordinates, we see that each four-derivation
\begin{align*}
\eta:\left(\mathcal{C}^\infty(\mathbb{R}^n)\right)^4\to\mathcal{C}^\infty(\mathbb{R}^n)
\end{align*}
has a tensor representation
\begin{align*}
\eta(f,s,h,q) = \sum_{i,j,k,\ell = 1}^n \eta_{i,j,k,\ell}\mathrm{d}_i f\mathrm{d}_j s\mathrm{d}_k h\mathrm{d}_\ell q
\end{align*}
for all $f,s,h,q\in\mathcal{C}^\infty(\mathbb{R}^n)$. Thus, $e$ has a representation $\varepsilon\in\mathcal{C}^\infty(\mathbb{R}^n,\mathbb{R}^{n\times n\times n\times n})$ with~\eqref{eq:conservative_tensor}. It remains to prove that $e$ has the properties (ii)--(iv) from Definition~\ref{def:cons_irrev} if, and only if, $\varepsilon$ has the properties (a)--(c).

\textsc{Step 2:} We show that $e$ fullfills the symmetry condition (ii) if, and only if, its tensor representation $\varepsilon$ fulfills (a). Let $i,j,k,\ell\in\underline{n}$ and assume that $e$ fulfills (ii). By evaluating $e$ at the coordinate functions, we find
\begin{align*}
\varepsilon_{i,j,k,\ell} = e(x_i,x_j,x_k,x_\ell) = e(x_i,x_j,x_\ell,x_k) = \varepsilon_{i,j,\ell,k}.
\end{align*}
Thus, $\varepsilon$ fulfills (a). The converse implication can be equally easy proven.

\textsc{Step 3:} We show that $e$ fulfills (iii) if, and only if,
\begin{align}\label{eq:bedingung_1}
\varepsilon_{i,j,i,\ell}+\varepsilon_{i,j,\ell,i}+\varepsilon_{\ell,j,i,i} = 0
\end{align}
for all $i,j,\ell\in\underline{n}$ and
\begin{align}\label{eq:Bedingung_2}
\varepsilon_{i,j,k,\ell} + \varepsilon_{k,j,i,\ell}+\varepsilon_{k,j,\ell,i}+\varepsilon_{\ell,j,k,i}+\varepsilon_{i,j,\ell,k}+\varepsilon_{\ell,j,i,k} = 0
\end{align}
for all $j\in\underline{n}$ and pairwise different $jik,\ell\in\underline{n}$.

``$\implies$'' We show first necessity. By evaluating $e$ at the coordinate functions, we get
\begin{align*}
\forall i,\ell\in\underline{n}: 0 = e(x_\ell,x_i,x_\ell,x_\ell) = \varepsilon_{\ell,i,\ell,\ell}.
\end{align*}
Let now $b\in\mathbb{R}\setminus\{0\}$. We have
\begin{align*}
0 & = e(x_i+bx_j,x_k,x_i+bx_j,x_i+bx_j)\\
& = b \varepsilon_{i,k,i,j}+b\varepsilon_{i,k,j,i}+b\varepsilon_{j,k,i,i}+b^2\varepsilon_{j,k,j,i}\\
& \quad +b^2\varepsilon_{j,k,i,j}+b^2\varepsilon_{i,k,j,j}
\end{align*}
for all $i\neq j\in\underline{n}$ and $k\in\underline{n}$. If we divide by $b$ and let $b$ tend to zero, we get
\begin{align*}
\forall i,j,k\in\underline{n}: 0 = \varepsilon_{i,i,j}+\varepsilon_{i,j,i}+\varepsilon_{j,i,i}
\end{align*}
and hence~\eqref{eq:bedingung_1} is indeed necessary. The property~\eqref{eq:Bedingung_2} can now be concluded by calculating $e(x_i+x_k+x_\ell,x_j,x_i+x_k+x_\ell,x_i+x_k+x_\ell)$ for $j\in\underline{n}$ and pairwise different $i,j,k\in\underline{n}$; we omit the details.

``$\impliedby$'' It is now straightforward to show sufficiency. First, we note that~\eqref{eq:bedingung_1} for $i = \ell$ yields $\varepsilon_{i,j,i,i} = 0$ for all $i,j\in\underline{n}$. Let $s,h\in\mathcal{C}^\infty(\mathbb{R}^n)$ be arbitrary. Then we have
\begin{align*}
e(h,s,h,h) & = \sum_{i,j,k,\ell = 1}^n \varepsilon_{i,j,k,\ell}\mathrm{d}_i h \mathrm{d}_j s \mathrm{d}_k h \mathrm{d}_\ell h\\
& = \sum_{j = 1}^n \mathrm{d}_j s\Bigg(\sum_{i = 1}^n\varepsilon_{i,j,i,i}(\mathrm{d}_i h)^3\\
& \quad + \sum_{i\neq k = 1}^n(\varepsilon_{i,j,i,k}+\varepsilon_{i,j,k,i}+\varepsilon_{k,j,i,i}) (\mathrm{d}_i h)^2\mathrm{d}_j h\\
&\quad + \sum_{i = 1}^n\sum_{k = i+1}^n\sum_{\ell = k+1}^n \varepsilon_{i,k,\ell} \mathrm{d}_i h \mathrm{d}_k h\mathrm{d}_\ell h\Bigg)\\
& = \sum_{j = 1}^n\mathrm{d}_j s\sum_{i = 1}^n\sum_{k = i+1}^n\sum_{\ell = k+1}^n (\varepsilon_{i,j,k,\ell}+\varepsilon_{i,j,\ell,k}\\
& \quad+\varepsilon_{k,j,i,\ell}+\varepsilon_{k,j,\ell,i}+\varepsilon_{\ell,j,i,k}+\varepsilon_{\ell,j,k,i}) \mathrm{d}_i h \mathrm{d}_j h\mathrm{d}_k h\\
& = 0.
\end{align*}
This shows the assertion.

\textsc{Step 4:} We show that, given the symmetry condition (a),~\eqref{eq:bedingung_1} and~\eqref{eq:Bedingung_2} are equivalent to (b). When we plug the symmetry condition (a) in, then we see that~\eqref{eq:Bedingung_2} is equivalent to
\begin{align*}
\varepsilon_{i,j,k,\ell}+\varepsilon_{k,j,\ell,i}+\varepsilon_{\ell,j,i,k} = 0
\end{align*}
for all $j\in\underline{n}$ and pairwise different $i,k,\ell\in\underline{n}$. Combining this with~\eqref{eq:bedingung_1}, we see that~\eqref{eq:bedingung_1} and ~\eqref{eq:Bedingung_2} are, in view of (a), indeed equivalent to (b).

\textsc{Step 5:} From the relation~\eqref{eq:conservative_tensor}, it is immediately clear that (iv) and (c) are equivalent. This shows that the proposition holds indeed true.\hfill$\square$

In our following considerations, we identify therefore the four-derivations $e$ from definition~\ref{def:cons_irrev} with their tensorial representations $\varepsilon\in\mathcal{C}^\infty(\mathbb{R}^{n},\mathbb{R}^{n\times n\times n\times n})$.

\section{Relation to port-Hamiltonian systems}

We have introduced conservative-irreversible functions. In Example~\ref{ex:an_example}, we have seen that IPHS can be geometrically formulated as being generated by conservative-irreversible functions according to \eqref{AutonomousIPHS}. In this section, we shall be interested in the converse question: which are the  conservative-irreversible functions which may be related to the expression in \eqref{AutonomousIPHS} generated by a quasi-Poisson bracket, i.e. Poisson bracket not satisfying the Jacobi identity? 

Before investigating this problem, let us consider the two-dimensional case and give an example of a conservative-irreversible function from which one may derive such a bracket.

\begin{Example}\label{ex:cons-irrev_ist_mehr}
Let $n = 2$ and define $\varepsilon\in\mathbb{R}^{2\times 2\times 2\times 2}$ by
\begin{align*}
\varepsilon_{i,j,k,\ell} := \begin{cases}
2, & (i,j,k,\ell)\in\{(1,1,2,2),(2,2,1,1)\},\\
-1, & (i,j,k,\ell)\in\left\{\begin{array}{l}
(1,2,1,2),(1,2,2,1),\\
(2,1,1,2),(2,1,2,1)
\end{array}\right\},\\
0, & \text{else}.
\end{cases}
\end{align*}
$\varepsilon$ is symmetric in the first two arguments and in the last two arguments. It fullfills, by construction, the condition (b) in Proposition~\ref{prop:tensorial_representation}. When checking condition (c) in Proposition~\ref{prop:tensorial_representation}, we have to consider, for all $h\in\mathcal{C}^\infty(\mathbb{R}^2)$,
\begin{align*}
& \begin{bmatrix}
\sum_{k,\ell = 1}^2\varepsilon_{1,1,k,\ell}\mathrm{d}_kh\mathrm{d}_\ell h & \sum_{k,\ell = 1}^2\varepsilon_{1,2,k,\ell}\mathrm{d}_kh\mathrm{d}_\ell h\\
\sum_{k,\ell = 1}^2\varepsilon_{2,1,k,\ell}\mathrm{d}_kh\mathrm{d}_\ell h & \sum_{k,\ell = 1}^2\varepsilon_{2,2,k,\ell}\mathrm{d}_kh\mathrm{d}_\ell h
\end{bmatrix}\\
 & = \begin{bmatrix}
2 \mathrm{d}_2h^2 & -2\mathrm{d}_1h\mathrm{d}_2h\\
-2\mathrm{d}_1h\mathrm{d}_2h & 2\mathrm{d}_1h^2
\end{bmatrix} =: M(h).
\end{align*}
We can calculate the characteristic polynomial of $M(h)$ as
\begin{align*}
p_{M(x)}(\lambda) & = (\lambda-2\mathrm{d}_2h^2)(\lambda-2\mathrm{d}_1h^2)-4\mathrm{d}_1h^2\mathrm{d}_2h^2\\
& = \lambda^2-2\lambda(\mathrm{d}_2h^2+\mathrm{d}_1h^2)
\end{align*}
and hence we conclude that $M(h)$ has eigenvalues $0$ and $2\mathrm{d}_1h^2+2\mathrm{d}_2h^2$. This shows that $M(h)$ is pointwise symmetric and positive semidefinite. Therefore, we conclude with Proposition~\ref{prop:tensorial_representation} that $(\varepsilon_{i,j,k,\ell})$ defines a conservative-irreversible function. \\
We relate this tensor to a quasi-Poisson bracket. The vector space of skew-symmetric matrices is one-dimensional and spanned by 
\begin{align*}
J := \begin{bmatrix}
0 & 1\\
-1 & 0
\end{bmatrix}.
\end{align*}
A tensorial representation $e\in\mathcal{C}^\infty(\mathbb{R}^n,\mathbb{R}^{n\times n\times n\times n})$ of the conservative-irreversible function
\begin{align*}
(f,s,h)\mapsto\{f,h\}_{J}\{s,h\}_J
\end{align*}
can thus be calculated as
\begin{align*}
\forall i,j,k,\ell\in\underline{2}: e_{i,j,k,\ell} = J_{i,k}J_{j,\ell}
\end{align*}
and hence we have
\begin{align*}
e_{i,j,k,\ell} = \begin{cases}
1, & (i,j,k,\ell)\in\{(1,1,2,2),(2,2,1,1)\},\\
-1, & (i,j,k,\ell) \in\{(1,2,2,1),(2,1,1,2)\},\\
0, & \text{else}.
\end{cases}
\end{align*}
When we symmetrise $2e$ in the last two entries, i.e. consider the four-tensor
\begin{align*}
[e_{i,j,k,\ell}+e_{i,j,\ell,k}]_{i,j,k,\ell = 1}^2,
\end{align*}
then we see that $\varepsilon$ is precisely this symmetrisation. Thus, the conservative-irreversible function represented by $\varepsilon$ is
\begin{align*}
(f,s,h)\mapsto 2\{s,h\}_J\{f,h\}_J.
\end{align*}
\hfill$\diamond$
\end{Example}

This example leads to study the set of conservative-irreversible functions that can be split in the form
\begin{align}\label{eq:Poisson_splitting}
(f,s,h)\mapsto\{s,h\}_1\{f,h\}_2
\end{align}
for biderivations $\{\cdot,\cdot\}_1,\{\cdot,\cdot\}_2$. This yields the question: Which conservative-irreversible functions allow a splitting~\eqref{eq:Poisson_splitting}?

\begin{Theorem}
Let $E:\left(\mathcal{C}^\infty(\mathbb{R}^n)\right)^3\to\mathcal{C}^\infty(\mathbb{R}^n)$ be a conservative-irreversible function. Then $E$ has the representation
\begin{align}\label{eq:bracket_splitting}
E(f,s,h) = \{s,h\}_1\{f,h\}_2
\end{align}
for some biderivations $\{\cdot,\cdot\}_1$, $\{\cdot,\cdot\}_2$ if, and only if, there is a quasi-Poisson bracket $\{\cdot,\cdot\}_J$ and $\gamma\in\mathcal{C}^\infty(\mathbb{R}^n,\mathbb{R}_{\geq 0})$ so that
\begin{align}\label{eq:pseudo_Poisson_splitting}
E(f,s,h) = \gamma\{s,h\}_J\{f,h\}_J.
\end{align}
\end{Theorem}
\textit{Proof.} Sufficiency of~\eqref{eq:pseudo_Poisson_splitting} follows from Example~\ref{ex:an_example}; we show necessity. Let $E$ have the representation~\eqref{eq:bracket_splitting} and the tensorial representation $e\in\mathcal{C}^\infty(\mathbb{R}^n,\mathbb{R}^{n\times n\times n\times n})$. The biderivations $\{\cdot,\cdot\}_1,\{\cdot,\cdot\}_2$ have the representation
\begin{align*}
\forall f,g\in\mathcal{C}^\infty(\mathbb{R}^n):\Bigg\{\begin{array}{ll}
\{f,g\}_1 & = \mathrm{d} f^\top A\mathrm{d} g,\\
\{f,g\}_2 & = \mathrm{d} f^\top B\mathrm{d} g
\end{array}
\end{align*}
for some $A,B\in\mathcal{C}^\infty(\mathbb{R}^n,\mathbb{R}^{n\times n})$. Then, $e$ can be calculated as
\begin{align}\label{eq:A_B_induce_e}
\forall i,j,k,\ell\in\underline{n}: e_{i,j,k,\ell} = A_{i,k}B_{j,\ell}.
\end{align}
Since the properties (a)--(c) in Proposition~\ref{prop:tensorial_representation} are in a pointwise manner, we can without loss of generality assume that $A$ and $B$ (and therefore $e$) are constant. Further, if either of the brackets is trivial, then the proposition holds obviously true, so that we can restrict ourself without loss of generality to nontrivial biderivations.

We show: Given matrices $A,B\in\mathbb{R}^{n\times n}\in\mathbb{R}^n\setminus\{0\}$, the four-tensor $e := [A_{i,k}B_{j,\ell}]_{i,j,k,\ell = 1}^n\in\mathbb{R}^{n\times n\times n\times n}$ fulfills (b) and (c) -- equivalently (iii) and (iv) -- if, and only if, $A$ is skew-symmetric and $A$ and $B$ are linearly dependent. We split the proof into steps.

\textsc{Step 1:} In the proof of Proposition~\ref{prop:tensorial_representation}, we have shown that the property (iii) of $e$ is equivalent to~\eqref{eq:bedingung_1} and~\eqref{eq:Bedingung_2}. We prove that
\begin{align}\label{eq:Bedingung_1_neu}
\forall i,j,k\in\underline{n}: e_{i,j,i,k}+e_{i,j,k,i}+e_{k,j,i,i} = 0\tag{7}
\end{align}
and
\begin{equation}\label{eq:Bedingung_2_neu}
\begin{aligned}
\forall j\in\underline{n}~\forall i,k,\ell\in\underline{n}~\text{pairwise~different}: e_{i,j,k,\ell}+e_{k,j,i,\ell}\\+e_{i,j,\ell,k}+e_{\ell,j,i,k}+e_{k,j,\ell,i}+e_{\ell,j,k,i} = 0
\end{aligned}\tag{8}
\end{equation}
hold true if, and only if, $A$ is skew-symmetric.

``$\implies$'' If we plug the definition of $e_{i,j,k,\ell}$ into~\eqref{eq:Bedingung_1_neu}, then we have
\begin{align*}
\forall i,j,k\in\underline{n}: A_{i,i}B_{j,k}+B_{j,i}(A_{i,k}+A_{k,i}) = 0.
\end{align*}
Especially, for $i = k$, we get
\begin{align*}
\forall i,j\in\underline{n}: A_{i,i}B_{j,i} = 0,
\end{align*}
so that we conclude
\begin{align*}
\forall i\in\underline{n}: A_{i,i} = 0\vee B_{\cdot,i} = 0.
\end{align*}
If $B_{\cdot,i} = 0$ for some fixed $i\in\underline{n}$, then we get 
\begin{align*}
\forall j,k\in\underline{n}: A_{i,i}B_{j,k} = 0
\end{align*}
and hence, since we assume $B\neq 0$, $A_{i,i} = 0$ in any case. Thus, we derive
\begin{align*}
\forall i,j,k\in\underline{n}: B_{j,i}(A_{i,k}+A_{k,i}) = 0
\end{align*}
from which we derive that $B_{\cdot,i} = 0$ or $A_{i,\cdot} = -A_{\cdot,i}$ for all $i\in\underline{n}$. If we plug this find into (ii), then we get, for all $j\in\underline{n}$ and all pairwise different $i,k,\ell\in\underline{n}$,
\begin{align*}
0 & = (A_{i,k}+A_{k,i})B_{j,\ell}+(A_{i,\ell}+A_{\ell,i})B_{j,k}+(A_{k,\ell}+A_{\ell,k})B_{j,i}.
\end{align*}
If $A$ is not skew-symmetric, then we have some $i\neq \ell\in\underline{n}$ so that $A_{i,\ell}\neq - A_{\ell,i}$ and thus $B_{j,i} = B_{j,\ell} = 0$ for all $j\in\underline{n}$. Hence we conclude 
\begin{align*}
\forall j\in\underline{n}~\forall k\in\underline{n}\setminus\{i,j\}: 0 = (A_{i,\ell}+A_{\ell,i})B_{j,k}.
\end{align*}
Therefore
\begin{align*}
\forall j\in\underline{n}~\forall k\in\underline{n}\setminus\{i,j\}: B_{j,k} = 0,
\end{align*}
so that we conclude with our previous findings that $B = 0$, contrary to our assumption. Hence $A$ is indeed skew-symmetric.

``$\impliedby$'' This follows from a straightforward calculation using the definition of $e$. We have, for all $i,j,k,\ell\in\underline{n}$
\begin{align*}
& e_{i,j,k,\ell}+e_{k,j,i,\ell}+e_{i,j,\ell,k}+e_{\ell,j,i,k}+e_{k,j,\ell,i}+e_{\ell,j,k,i}\\
&\quad = (A_{i,k}+A_{k,i})B_{j,\ell}+(A_{i,\ell}+A_{\ell,i})B_{j,k}+(A_{k,\ell}+A_{\ell,k})B_{j,i}\\
&\quad = 0
\end{align*}
by skew-symmetry of $A$ and hence~\eqref{eq:Bedingung_2_neu} holds true. Putting $i = k$, we readily conclude that~\eqref{eq:Bedingung_1_neu} holds equally true.

This proves that $e$ has the property (iii) if, and only if, $A$ is skew-symmetric.

\textsc{Step 2:} It is easy to verify that the symmetry of $e(\cdot,\cdot,h,h)$ in property (iii) is equivalent to
\begin{align}\label{eq:symmetry_property}
\forall y\in\mathbb{R}^n: Ayy^\top B^\top = Byy^\top A^\top.
\end{align}
We show: If $A$ and $B$ fulfill~\eqref{eq:symmetry_property},
then $A_{\cdot,i}$ and $B_{\cdot,i}$ are linearly dependent for all $i\in\underline{n}$. Especially, if $A$ and $B$ fulfill~\eqref{eq:symmetry_property}, then, for each $i\in\underline{n}$, the matrix
\begin{align*}
Ae_ie_i^\top B^\top = A_{\cdot,i}B_{\cdot,i}^\top,
\end{align*}
where $e_i$ denotes, as usual, the $i$-th standard unit vector, is symmetric. Therefore, we conclude
\begin{align*}
\forall\lambda\in\mathbb{R}~\forall j,k\in\underline{n}: B_{j,i}(A_{k,i}-\lambda B_{k,i}) = B_{i,k}(A_{j,i}-\lambda B_{j,i});
\end{align*}
the same holds true with $A$ and $B$ interchanged. If $B_{\cdot,i}\neq 0$, then there is some $k\in\underline{n}$ with $B_{k,i}\neq 0$ and some $\mu_i\in\mathbb{R}$ so that $A_{k,i} = \mu_i B_{k,i}$ and thus we have
\begin{align*}
\forall j\in\underline{n}: 0 = B_{j,i}(A_{k,i}-\mu_i B_{k,i}) = B_{k,i}(A_{j,i}-\mu_i B_{j,i}),
\end{align*}
so that we conclude $A_{\cdot,i} = \mu_i B_{\cdot,i}$. Interchanging $A$ and $B$, we conclude that if $A_{\cdot,i}\neq 0$, then there is some $\lambda_i\in\mathbb{R}$ so that $B_{\cdot,i} = \lambda_i A_i$. If, on the other hand, $B_{\cdot,i} = 0$ and $A_{\cdot,i} = 0$, then $A_{\cdot,i}$ and $B_{\cdot,i}$ are linearly dependent, anyway. Therefore,
\begin{align*}
\forall i\in\underline{n}~\exists (\lambda_i,\mu_i)\in\mathbb{R}^2\setminus\{0\}: \lambda_i A_{\cdot,i} + \mu_i B_{\cdot,i} = 0.
\end{align*}

\textsc{Step 3:} In this step, we exploit the skew-symmetry of $A$, which we verified in \textsc{Step 1}. We show that if $A$ and $B$ fullfill~\eqref{eq:symmetry_property} and $A\neq 0$ is skew-symmetric, then $Ayy^\top B^\top$ is symmetric for all $y\in\mathbb{R}^n$ if, and only if, $A$ and $B$ are linerly dependent. Let $i,j\in\underline{n}$ be so that $A_{i,j}\neq 0$ (especially $i\neq j$); after renormalisation we can w.l.o.g. assume that $A_{i,j} = 1$. Then $A_{i,j} = -A_{j,i}$ so that we can w.l.o.g. further assume that $i<j$ and thus
\begin{align*}
A(e_i+e_j) = 
\left(\begin{array}{ccccccccccc}
\star & \cdots & \star & 1 & \star & \cdots & \star & -1 & \star & \cdots & \star
\end{array}\right)^\top
\end{align*}
By \textsc{Step 2} there are $\lambda_i,\lambda_j\in\mathbb{R}$ so that $B_{\cdot,i} = \lambda_i A_{\cdot,i}$ and $B_{\cdot,j} = \lambda_j A_{\cdot,j}$. Plugging this in, we get
\begin{align*}
\big(A(e_i+e_j)(e_i+e_j)^\top B^\top\big)_{i,j} = \lambda_i
\end{align*}
and
\begin{align*}
\big(A(\lambda e_i+\mu e_j)(\lambda e_i+\mu e_j)^\top B^\top\big)_{j,i} = \lambda_j.
\end{align*}
By~\eqref{eq:symmetry_property}, we conclude $\lambda_i = \lambda_j$. Define the graph $G = (\underline{n},E)$ with
\begin{align*}
(\iota,\kappa)\in E :\iff \big(\exists i\in\underline{n}: A_{\iota
,i},A_{\kappa,i}\neq 0\big)\vee A_{\iota,\kappa}\neq 0
\end{align*}
and the equivalence relation $\sim$ that is canonically associated to $G$ and induced by the connected components of $G$. Then we conclude for all $i,j\in\underline{n}$ with $i\in[j]_{/\sim}$ that $\lambda_i = \lambda_j$. It remains to consider the case $A_{\cdot,i} = 0$. In this case, we have
\begin{align*}
A(e_i+e_j)(e_i+e_j)^\top B^\top = Ae_ie_i^\top B^\top + A_{\cdot,j}B_{\cdot,i}^\top.
\end{align*}
Thus, $A_{\cdot,j}B_{\cdot,i}^\top$ is symmetric. Analogously to \textsc{Step 2}, we conclude that for all $j\in\underline{n}$ with $A_{\cdot,j} \neq 0$, there is some $\lambda_{i,j}\in\mathbb{R}$ so that $B_{\cdot,i} = \lambda_{i,j}A_{\cdot,j}$. Since $A\neq 0$, and since $A$ is skew-symmetric, $A$ contains at least two (non-zero) linearly independent columns. Thus, we have
\begin{align*}
B_{\cdot,i}\in\bigcap_{i\in\underline{n}, A_{\cdot,i}\neq 0}\mathrm{span}\{A_{\cdot,i}\} = \{0\}.
\end{align*}
Hence, after possibly performing some simultaneous permutations of the rows and likewise the columns of $A$ and $B$ -- which are orthogonal transformations that do not change the symmetry properties of $A$ and $B$ -- we can split
\begin{align*}
A = \begin{bmatrix}
A^1\\
& A^2\\
& & \ddots\\
& & & A^k
\end{bmatrix},\qquad B = \begin{bmatrix}
B^1\\
& B^2\\
& & \ddots\\
& & & B^k
\end{bmatrix}
\end{align*}
for some $k\in\underline{n}$, where the $A^i\in\mathbb{R}^{n_i\times n_i}$ and $B^i\in\mathbb{R}^{n_i\times n_i}$ have the form
\begin{align*}
A^i & = \begin{bmatrix}
& & & & & & A^i_1\\
& & & & & \udots\\
& & & & A^i_{\ell_i}\\
& & & 0\\
& & -(A^i_{\ell_i})^\top\\
& \udots\\
-(A^i_1)^\top
\end{bmatrix},\\
B^i & = \begin{bmatrix}
& & & & & & \lambda_{i,1}A^i_1\\
& & & & & \udots\\
& & & & \lambda_{i,\ell_i}A^i_{\ell_i}\\
& & & 0\\
& & -\lambda_{i,\ell_i}(A^i_{\ell_i})^\top\\
& \udots\\
-\lambda_{i,1}(A^i_1)^\top
\end{bmatrix}
\end{align*}
We show first that $\lambda_{i,j} = \lambda_{i,k}$ for all $i\in\underline{n}$ and $j,k\in\underline{\ell_i}$. Let $y\in\mathbb{R}^{n_i}$ be arbitrary. Then we have
\begin{align*}
A^iyy^\top B^i & = \begin{pmatrix}
A^i_1y_n\\\vdots\\A^i_{\ell_i}y_{\ell_i + 2}\\0\\-A^i_{\ell_i}y_{\ell_i}\\\vdots\\-A^i_1y_1
\end{pmatrix}\begin{pmatrix}
\lambda_{i,j}y_{n-j+1}^\top(A^i_1)^\top
\end{pmatrix}_{j = 1}^n
\end{align*}
which, in turn, equals the rather lengthy matrix depictured in Figure~\ref{eq:stupid_equation_2}.
\begin{figure*}\label{eq:stupid_equation_2}
\begin{align}\label{eq:stupid_equation1}
\begin{bmatrix}
\lambda_{i,1} A^i_1y_ny_n^\top(A^i_1)^\top & \cdots & \lambda_{i,\ell_1}A^i_1y_ny_{\ell_i + 2}^\top (A^i_{\ell_i})^\top & 0 & -\lambda_{i,\ell_1}A^i_1y_ny_{\ell_i}^\top (A^i_{\ell_i})^\top & \cdots & -\lambda_{i,1}A^i_1y_ny_1^\top (A^i_1)^\top\\
\vdots & \ddots & \vdots & \vdots & \vdots & \ddots & \vdots\\
\lambda_{i,1} A^i_{\ell_i}y_{\ell_i + 2}y_n^\top(A^i_1)^\top & \cdots & \lambda_{i,\ell_1}A^i_{\ell_i}y_{\ell_i + 2}y_{\ell_i + 2}^\top (A^i_{\ell_i})^\top & 0 & -\lambda_{i,\ell_1}A^i_{\ell_i}y_{\ell_i + 2}y_{\ell_i}^\top (A^i_{\ell_i})^\top & \cdots & -\lambda_{i,1}A^i_{\ell_i}y_{\ell_i + 2}y_1^\top (A^i_1)^\top\\
0 & \cdots & 0 & 0 & 0 & \cdots & 0\\
-\lambda_{i,1} A^i_{\ell_i}y_{\ell_i}y_n^\top(A^i_1)^\top & \cdots & -\lambda_{i,\ell_1}A^i_{\ell_i}y_{\ell_i}y_{\ell_i + 2}^\top (A^i_{\ell_i})^\top & 0 & \lambda_{i,\ell_1}A^i_{\ell_i}y_{\ell_i}y_{\ell_i}^\top (A^i_{\ell_i})^\top & \cdots & \lambda_{i,1}A^i_{\ell_i}y_{\ell_i}y_1^\top (A^i_1)^\top\\
\vdots & \ddots & \vdots & \vdots & \vdots & \ddots & \vdots\\
-\lambda_{i,1} A^i_1y_1y_n^\top(A^i_1)^\top & \cdots & -\lambda_{i,\ell_1}A^i_1y_1y_{\ell_i + 2}^\top (A^i_{\ell_i})^\top & 0 & \lambda_{i,\ell_1}A^i_1y_1y_{\ell_i}^\top (A^i_{\ell_i})^\top & \cdots & \lambda_{i,1}A^i_1y_1y_1^\top (A^i_1)^\top
\end{bmatrix}.
\end{align}
\caption{Structure of $A^iyy^\top B^i$.}
\end{figure*}
Since this matrix is symmetric by assumption, we conclude that indeed for all $j,k\in\underline{\ell_i}$ either $\lambda_{i,j} = \lambda_{i,k}$ or $A^i_jy_{n-j+1}y_{n-k+1}^\top (A^i_k)^\top = 0$ for all $y_{n-j+1},y_{n-k+1}$; but if the latter is the case, then $A^i_j = 0$ or $A^i_k = 0$ and we can w.l.o.g. assume that $\lambda_{i,j} = \lambda_{i,k}$ anyway. This shows that it remains to consider the case
\begin{align*}
A = \begin{bmatrix}
A^1\\
& \ddots\\
& & A^k
\end{bmatrix},\quad B = \begin{bmatrix}
\lambda_1 A^1\\
& \ddots\\
& & \lambda_k A^k
\end{bmatrix}
\end{align*}
for some $k\in\mathbb{N}^*$ and $\lambda_1,\ldots,\lambda_k\in\mathbb{R}$. If $k = 1$, then $A$ and $B$ are linearly dependent. Consider the case $k = 2$. Then we have for all $y = (y_1,y_2)\in\mathbb{R}^n$
\begin{align*}
Ayy^\top B^\top = \begin{bmatrix}
\lambda_1 A_1y_1y_1^\top A_1^\top & \lambda_2 A_1y_1y_2^\top A_2^\top\\[0.5\normalbaselineskip]
\lambda_1 A_2y_2y_1^\top A_1^\top & \lambda_2 A_2y_2y_2^\top A_2^\top
\end{bmatrix}
\end{align*}
and hence we conclude $(\lambda_1-\lambda_2)A_2y_2y_1^\top A_1^\top = 0$. This yields $A_1 = 0$ or $A_2 = 0$ or $\lambda_1 = \lambda_2$. But if $A_1 = 0$ or $A_2 = 0$, then we can choose $\lambda_1 = \lambda_2$ anyways, so that we conclude that w.l.o.g. $\lambda_1 = \lambda_2$. By induction over $k$, we conclude that indeed $\lambda_1 = \cdots = \lambda_k$. This concludes the proof of the theorem.\hfill $\square$

We have seen that the conservative-irreversible functions induced by quasi-IPHS are all conservative-irreversible functions that allow a splitting~\eqref{eq:bracket_splitting}. Immediately, this raises the question: Are these functions $E\in\mathfrak{CI}(\mathbb{R}^n)$ all so that their associated four-derivation $e$ induces the family $(e(\cdot,s,h,\cdot))_{s,h\in\mathcal{C}^\infty(\mathbb{R}^n)}$ of quasi-Poisson brackets? Unfortunately, we have not the room to answer this question. We can, however, characterize the functions $E\in\mathfrak{CI}(\mathbb{R}^n)$ with this property.

\begin{prop}\label{prop:quasi_Hamiltonian_condition}
Let $E:\left(\mathcal{C}^\infty(\mathbb{R}^n)\right)^3\to\mathcal{C}^\infty(\mathbb{R}^n)$ be a conservative-irreversible function with associated four-derivation $e$. All biderivations $e(\cdot,s,h,\cdot)$ are quasi-Poisson if, and only if,
\begin{align}\label{eq:port_Hamiltonian_condition}
\forall i,j,k,\ell\in\underline{n}: e_{i,j,k,\ell} = -e_{\ell,j,k,i}.
\end{align}
\end{prop}
\textit{Proof.}
In view of Proposition~\ref{prop:tensorial_representation}, we identify $e$ with its representation $e\in\mathcal{C}^\infty(\mathbb{R}^n,\mathbb{R}^{n\times n\times n\times n})$. Then, it is evident that, for all $s,h\in\mathcal{C}^\infty(\mathbb{R}^n)$, the matrix representation of $e(\cdot,s,h,\cdot)$ is
\begin{align*}
[e(\cdot,s,h,\cdot)] = \sum_{j,k = 1}^ne_{i,j,k,\ell}\mathrm{d}_j s\mathrm{d}_k h.
\end{align*}
Thus, $e(\cdot,s,h,\cdot)$ is a quasi-Poisson bracket if, and only if,
\begin{align*}
\sum_{j,k = 1}^ne_{i,j,k,\ell}\mathrm{d}_j s\mathrm{d}_k h = -\sum_{j,k = 1}^ne_{\ell,j,k,i}\mathrm{d}_j s\mathrm{d}_k h.
\end{align*}
We show that latter and~\eqref{eq:port_Hamiltonian_condition} are equivalent. It is clear that~\eqref{eq:port_Hamiltonian_condition} is sufficient. We show necessity. By plugging in the coordinate functions $x_\iota$ and $x_\kappa$, we see that
\begin{align*}
e_{i,\iota,\kappa,\ell} & = 
\sum_{j,k = 1}^ne_{i,j,k,\ell}\mathrm{d}_j x_\iota\mathrm{d}_k x_\kappa\\
& = -\sum_{j,k = 1}^ne_{\ell,j,k,i}\mathrm{d}_j x_\iota\mathrm{d}_k x_\kappa\\
& = -e_{\ell,\iota,\kappa,i}.
\end{align*}
This shows that~\eqref{eq:port_Hamiltonian_condition} is indeed necessary.\hfill$\square$

\section{Conclusion and outlook}

We have suggested a geometric characterization of the quasi-Poisson brackets of Irreversible port-Hamiltonian Systems in terms of  covariant 4-tensors and an associated functions, called \emph{conservative-irreversible functions}. We have shown that all conservative-irreversible functions that can be split into two biderivations are induced by an quasi-Poisson bracket as defined for Irreversible port-Hamiltonian Systems. Lastly, we have characterized the conservative-irreversible functions that induce quasi-pseudo-Hamiltonian systems.

Unfortunately, it remains unanswered whether the conser\-vative-irreversible functions are precisely the functions that allow a splitting~\eqref{eq:Poisson_splitting}. In future works, we want to answer this question. Next, we have not taken account that the function $\gamma$ defining the dissipative constitutive relations for Irreversible port-Hamiltonian Systems, may depend explicitly on the co-energy variable. Note that, as long as this function is positive, the entropy creation term remain positive. Finally future work shall make use of these conservative-irreversible functions in order to define with an input- and an output-map  compatible with irreversibel thermodynamics. Finally, having characterized, possibly more general systems retaining both the energy and entropy balance equations, might lead to novel desired closed-loop systems useful for the control design.

\vspace*{1cm}
\begin{footnotesize}
\textbf{Acknowledgements} We thank our colleague Manuel Schaller (Ilmenau) for his valuable critique.
\end{footnotesize}

\bibliographystyle{alpha}
\bibliography{my_ifacconf} 
\end{document}